\documentclass[12pt]{amsart}
\usepackage{amssymb, amsmath}

\newtheorem{theorem}{Theorem}[section]

\newtheorem{lemma}[theorem]{Lemma}
\newtheorem*{theoremnn}{Theorem}

\theoremstyle{definition}

\theoremstyle{remark}

\begin{document}

\newcommand{\fra}[1]{{\mathfrak{#1}}}
\newcommand{\Gr}{\text{Gr}}
\newcommand{\grh}{{\text{gr}}_{F}}
\newcommand{\ra}{\rightarrow}
\newcommand{\Ra}{\Rightarrow}
\newcommand{\surj}{\twoheadrightarrow}
\newcommand{\lra}{\longrightarrow}
\newcommand{\noi}{\noindent}
\newcommand{\PP}{\mathbf{P}}
\newcommand{\PPP}{\PP_{\text{sub}}}
\newcommand{\RR}{\mathbf{R}}
\newcommand{\NN}{\mathbf{N}}
\newcommand{\lef}{\mathbf{L}}
\newcommand{\hyp}{\mathbf{H}}
\newcommand{\ZZ}{\mathbf{Z}}
\newcommand{\CC}{\mathbf{C}}
\newcommand{\QQ}{\mathbf{Q}}
\newcommand{\bin}[2]{ {{#1} \choose {#2}} }
\newcommand{\OO}{\mathcal{O}}

\newcommand{\cH}{\mathcal{H}}
\newcommand{\MM}{\mathcal{M}}
\newcommand{\KK}{\mathcal{K}}
\newcommand{\II}{\mathcal{I}}
\newcommand{\LL}{\mathcal{L}}
\newcommand{\FF}{\mathcal{F}}
\newcommand{\GG}{\mathcal{G}}
\newcommand{\EE}{\mathcal{E}}
\newcommand{\cM}{\mathcal{M}}
\newcommand{\cN}{\mathcal{N}}
\newcommand{\DD}{\mathcal{D}}
\newcommand{\cA}{\mathbf{A}}
\newcommand{\ti}[1]{\tilde{#1}}
\newcommand{\ef}{\rm{\ if\  }}
\newcommand{\eps}{\epsilon}
\newcommand{\sbl}{\vskip 3pt}
\newcommand{\lbl}{\vskip 6pt}
\newcommand{\rk} {\text{rank }}
\newcommand{\osect}{\mathbf{0}}
\newcommand{\HH}[3]{H^{{#1}} \big( {#2} , {#3} \big) }
\newcommand{\hh}[3]{h^{{#1}} \big( {#2} , {#3} \big) }
\newcommand{\OP}[1]{\OO_{\PP({#1})}}
\newcommand{\fall}{ \ \ \text{ for all } \ }
\newcommand{\rndup}[1]{ \ulcorner {#1} \urcorner }
\newcommand{\rndown}[1] {\llcorner {#1} \lrcorner}
\newcommand{\hgt}{\rm{height }}
\newcommand{\MI}[1]{\mathcal{J} ( {#1} ) }
\newcommand{\MIP}[1]{\mathcal{J}_+ ( {#1} ) }
\newcommand{\ZMI}[1]{\text{Zeroes}\big( \MI{ {#1} }
          \big)}
\newcommand{\BI}[1]{  \mathfrak{b} \big( {#1} \big) }
\newcommand{\ord}{\text{ord}}
\newcommand{\codim}{\text{codim}}
\newcommand{\mult}{\text{mult}}
\newcommand{\Supp}{\text{Supp}}
\newcommand{\defi}{{\text{def}}}
\newcommand{\pr}{\prime}
\newcommand{\QT}[2]{{#1}<{#2}>}
\newcommand{\Div}{\text{Div}}
\newcommand{\num}{ \equiv_{\text{num}} }
\newcommand{\lin}{\equiv}
\newcommand{\Qnum}{\equiv_{\text{num},\QQ}}
\newcommand{\Qlin}{\equiv_{\text{lin},\QQ}}
\newcommand{\dra}{\dashrightarrow}
\newcommand{\Bl}{\text{Bl}}
\newcommand{\lcm}{\text{lcm}}
\newcommand{\length}{\text{length}}
\newcommand{\lam}{\lambda}
\newcommand{\al}{\alpha}
\newcommand{\te}{\text}
\newcommand{\Spec}{\text{Spec}}
\newcommand{\HSp}{\text{HSp}}
\newcommand{\db}{\underline{\Omega}\,\dot{}}

\title[Multiplier ideals and filtered D-modules]{Multiplier ideals and filtered D-modules}
\author{Nero Budur}
\address{Department of Mathematics \\University of
Illinois at Chicago \hfil\break\indent 851 South Morgan Street
(M/C 249)\\ Chicago, IL 60607-7045, USA} \email{nero@math.uic.edu}
\date{May 15, 2003}

\begin{abstract}
We give a Hodge-theoretic interpretation of the multiplier ideal
of an effective divisor on a smooth complex variety. More
precisely, we show that the associated graded coherent sheaf with
respect to the jumping-number filtration can be recovered from the
smallest piece of M. Saito's Hodge filtration of the $\DD$-module
of vanishing cycles.
\end{abstract}

\maketitle

\section{Introduction}
\label{intro}

Let $D$ be an effective $\QQ$-divisor on a nonsingular complex
variety $X$ of dimension $n$. The  multiplier ideal $\MI{D}$ is a
subsheaf of ideals of $\OO_X$ and measures in a subtle way the
singularities of $D$, see \cite{La}. The singularities of $D$ get
"worse" if $\MI{D}$ is smaller. The main goal of this note is to
give a Hodge-theoretic interpretation of multiplier ideals. That
such an interpretation is possible was hinted by \cite{Bu} where
we proved a local relation at a point $x\in X$ between $\MI{D}$
and the mixed Hodge structure on the cohomology the Milnor fiber
of an integral divisor $D$ at $x$.

The natural setting for our result is the theory of mixed Hodge
modules due to M. Saito (\cite{Sa1}, \cite{Sa2}). Since we restrict our attention to the Hodge filtration only and disregard the weight filtration and the rational structure, we end up working with filtered $\DD_X$-modules $(M,F)$. Here $\DD_X$ is the sheaf of non-commutative rings of linear algebraic differential operators (see \cite{Bo}). The Hodge filtration $F$ is always assumed here to be increasing. By the Riemann-Hilbert correspondence, $M$ corresponds to a perverse sheaf on $X$. For example, the trivial mixed Hodge module $\QQ_X^H[n]$ is represented by the filtered left $\DD_X$-module $(\OO_X,F)$, where $\Gr^F_p=0$ for $p\ne 0$. The corresponding perverse sheaf is the shifted trivial complex $\QQ_X[n]$.

For a non-constant regular function $f:X\ra\CC$, one has the vanishing cycles functor $\psi_f$ which can be defined on the abelian category of mixed Hodge modules. By definition, $\psi_f$ corresponds to
${}^p\psi_f=\psi_f[-1]$ on the category of perverse sheaves on
$X$. For $\al\in (0,1]\cap \QQ$, let $\psi_{f}^\al\OO_X$
correspond to the eigenspace of the semisimple part of monodromy
for the eigenvalue $\exp(-2\pi i\al)$. M. Saito's theory provides
us with a canonical filtration $F$ on $\psi_{f}^\al\OO_X$ (for definitions, see the introduction of \cite{Sa1}).

\begin{theoremnn}\label{thm}
Let $X$ be a smooth complex variety of dimension $n$. Let
$f:X\ra\CC$ be a non-constant regular function and $D=f^{-1}(0)$
the corresponding effective divisor. Then for $\al\in (0,1]$,
\begin{align}\label{t1}
\frac{\MI{(\al-\eps)\cdot D}}{\MI{\al\cdot
D}}=F_{1}\psi_{f}^\al\OO_X,
\end{align}
where $0<\eps\ll 1$.
\end{theoremnn}

Here, $F_1$ is the smallest piece of the Hodge filtration of the left $\DD_X$-module $\psi_f\OO_X$. The values $\al\in (0,1]$ for which the left-hand side of
(\ref{t1}) is nonzero are called jumping numbers (see \cite{Bu},
\cite{La}, \cite{ELSV}). The  values $\al\in (0,1]$ for which the
right-hand side of (\ref{t1}) is nonzero were considered in
\cite{Sa93}. Thus the Theorem answers a question in \cite{ELSV}
regarding the relation between the two sets of values and reproves
their theorem that the jumping numbers of $D$ are roots of the
Bernstein-Sato polynomial of $f$ up to a sign.

I thank Professor L. Ein for several discussions and comments on this work.

{\it Acknowledgement.} Professor M. Saito has also proved the Theorem
of this note, see \cite{Sa03}. Moreover, he proves that the
multiplier ideals $\MI{\al\cdot D}$ give the $V$-filtration of
Malgrange and Kashiwara corresponding to $f$ on the left $\DD_X$-module
 $\OO_X$.

\section{Proof of the Theorem}

We used left $\DD$-modules only for the introduction. We will work, as our references do, with right $\DD_X$-modules. The trivial right $\DD_X$-module is $\omega_X=\bigwedge^n\Omega_X^1$, the sheaf of regular $n$-forms. Locally, the action of a vector field $\xi$ on $w\in\omega_X$ is given by $w\xi=-\te{Lie}_\xi w$,  the Lie derivative. $\QQ_X^H[n]$ is represented by the filtered right $\DD_X$-module $(\omega_X,F)$, where $\Gr^F_p=0$ for $p\ne -n$. In general, the transformation from left to right $\DD_X$-modules is given by
\begin{align*}
(M,F)&\mapsto (\ \omega_X\otimes M, \  F_p=\omega_X\otimes F_{p+n}\ ),\\
 &(w\otimes u)\xi=w\xi\otimes u-w\otimes\xi u.
\end{align*}
Hence (\ref{t1}) is equivalent to
\begin{align}\label{t}
\omega_X\otimes_{\OO_X}\frac{\MI{(\al-\eps)\cdot D}}{\MI{\al\cdot
D}}=F_{1-n}\psi_{f}^\al\omega_X.
\end{align}

 Let $\mu:Y\ra X$ be a log resolution of $(X,D)$. Recall that the multiplier ideal
 $\MI{\al\cdot D}$ is defined for all $\al>0$ by
$$\MI{\al\cdot D}=\mu_*\OO_Y(K_{Y/X}-\rndown{\al\cdot\mu^*D}),$$
where $K_{Y/X}=K_Y-\mu^*K_X$ and $\rndown{.}$ means rounding down
the coefficients in a $\QQ$-divisor. Put
$$\KK_\al(X,D)=\frac{\MI{(\al-\eps)\cdot D}}{\MI{\al\cdot D}},$$
where $0<\eps\ll 1$. Then $\MI{\al\cdot D}$ and $\KK_\al(X,D)$ are
independent of the choice of $\mu$ and
\begin{align}\label{kk}
\KK_\al(X,D)=\mu_*\left
(\OO_Y(K_{Y/X})\otimes\KK_\al(Y,\mu^*D)\right ).
\end{align}

Put $g=f\circ \mu$. Let $\MM=(\omega_Y,F)$ be the filtered
$\DD_Y$-module with $Gr_p^F=0$ for $p\ne -n$. By Theorem 2.14 of \cite{Sa2}, for $\al\in
(0,1]$,
$$\psi_f^\al H^j\mu_*\MM=H^j\mu_*\psi_g^\al\MM,$$
for all $j\in\ZZ$. Here $\mu_*:D^bMHM(Y)\ra D^bMHM(X)$ is the
direct image functor on the bounded derived categories of mixed
Hodge modules (we care only about the complexes of filtered
$\DD$-modules), and $H^j$ is the $j$-th cohomology of a complex.
In particular, we have an equality of filtered $\DD_X$-modules
\begin{align}\label{t2}
\psi_f^\al H^0\mu_*\MM=H^0\mu_*\psi_g^\al\MM.
\end{align}
We will show that (\ref{t}), hence the Theorem, follows by taking $F_{1-n}$ of both sides
of (\ref{t2}).

\begin{lemma}\label{l1}
$F_{1-n}\psi_f^\al H^0\mu_*\MM=F_{1-n}\psi_f^\al(\omega_X,F)$, for
$\al\in (0,1]$.
\end{lemma}
\begin{proof}
Follows from $H^0\mu_*\QQ_Y^H[n]=\QQ_X^H[n]$.
\end{proof}

\begin{lemma}\label{l2}
$F_{1-n}H^0\mu_*\psi_g^\al\MM=\omega_X\otimes_{\OO_X}\KK_\al(X,D),$
for $\al\in (0,1]$.
\end{lemma}
\begin{proof}
Recall the definition of $\mu_*(M,F)$. Let $i_\mu:Y\ra Y\times X$
be the graph of $\mu$. Let $p:Y\times X\ra X$ be the natural
projection. Then
$$\mu_*(M,F)=Rp_. DR_{X\times Y/X}(i_\mu)_*(M,F).$$
where $Rp_.$ is the usual derived direct image for sheaves. We put from now $p_.=H^0(Rp_.)$ for the usual direct image of sheaves. Recall
that $DR_{X\times Y/X}(M',F)$ is defined by
$$F_pDR_{X\times Y/X}(M',F)=\left
[F_{p-n}M'\otimes\bigwedge^n\Theta_Y\lra\ldots\lra F_pM'\right
],$$ where $F_pM'$ sits in degree zero in the last complex, and
$\Theta_Y=(\Omega_Y^1)^\vee$. The definition of
$(i_\mu)_*(M,F)=(M',F)$ is the same as for $\DD_Y$-modules, and
all we need to know about the Hodge filtration is that
$F_{q}M'=(i_\mu)_.F_{q}M$ if $q=\min\{\ p\ |\ F_pM\ne 0\ \}$. In
this case, also $q=\min\{\ p\ |\ F_pM'\ne 0\ \}$.

If $(M,F)=\psi_g^\al(\omega_Y,F)$ and $\al\in (0,1]$ is such that
$M\ne 0$, then $q=1-n$. Hence,
\begin{align*}
F_{1-n}H^0\mu_*\psi_g^\al (M,F)&=p_.(i_\mu)_.F_{1-n}M\\
&=\mu_.F_{1-n}M,
\end{align*}
where the last $\mu_.$ is the usual sheaf direct image. By Lemma
\ref{l3}, $F_{1-n}M=\omega_Y\otimes\KK_\al(Y,\mu^*D)$. By
(\ref{kk}), this proves the claim.
\end{proof}

\begin{lemma}\label{l3}
The Theorem is true if $D$ is a simple normal crossing divisor.
\end{lemma}
\begin{proof}
By definition, for $\al\in (0,1]$,
$$\psi_f^\al\omega_X=\Gr_V^\al(i_f)_*\omega_X,$$
where $i_f:X\ra X\times\CC$ is the graph of $f$, $(i_f)_*$ is the
direct image for (filtered) $\DD_X$-modules, and $V$ is the
decreasing filtration of Malgrange and Kashiwara. The Hodge
filtration on $\psi_f^\al$ is the filtration $F[1]$, where $F$ is
induced by $(i_f)_*$. Here $F[i]_p=F_{p-i}$. In particular,
$$F_{1-n}\psi_f^\al\omega_X=F_{-n}\Gr_V^\al(i_f)_*\omega_X.$$
By Proposition 3.5-(3.5.1) of \cite{Sa2} applied to $(\omega_X,
F[n])$, one has that
\begin{align*}
F_{-n}V^\al(i_f)_*\omega_X&=\omega_X\otimes\OO_X(-\rndown{(\al-\eps)\cdot
D})\\
&=\omega_X\otimes\MI{(\al-\eps)\cdot D},
\end{align*}
where $0<\eps\ll 1$. Indeed, to apply that Proposition one only
has to check locally, where $X$ has coordinates $x_1,\ldots, x_n$,
that $\omega(x_i\partial _i+1)=0$, for
$\omega=dx_1\wedge\ldots\wedge dx_n$, and $\partial
_i=\partial_{x_i}$. This follows from $\te{Lie}_{\partial
_i}(x_i\omega)=\omega$. Hence it gives $\mu=(-1,\ldots, -1)$ in the
above-mentioned Proposition.
\end{proof}

\end{document}